\documentclass[a4paper,12pt]{amsart}
\usepackage{amssymb}
\usepackage{ifthen}
\usepackage{graphicx}
\usepackage{color,xcolor}
\nonstopmode
\numberwithin{equation}{section}
\newtheorem{thm}{Theorem}
\newtheorem{cor}{Corollary}
\newtheorem{lem}{Lemma}
\newtheorem{prop}{Proposition}

\newtheorem{conj}{Conjecture}
\newtheorem{prob}{Problem}

\theoremstyle{definition}
\newtheorem{defn}{Definition}
\newtheorem{ca}{Case}
\newtheorem{example}{Example}
\newtheorem{rem}{Remark}

\newenvironment{pf}[1][]{%
 \vskip 1mm
 \noindent
 \ifthenelse{\equal{#1}{}}%
  {{\slshape Proof. }}%
  {{\slshape #1.} }%
 }%
{\qed\medskip}

\newcounter{alphabet}

\newenvironment{Thm}[1][]{\refstepcounter{alphabet}%
\bigskip%
\noindent%
{\bf Theorem \Alph{alphabet}}%
\ifthenelse{\equal{#1}{}}{}{ (#1)}%
{\bf .} \itshape}{\vskip 8pt}

\newenvironment{Lem}[1][]{\refstepcounter{alphabet}%
\bigskip%
\noindent%
{\bf Lemma \Alph{alphabet}}%
{\bf .} \itshape}{\vskip 8pt}

\newcounter{alphabet2}

\newcommand{\ID}{{\mathbb D}}

\newcommand{\real}{{\operatorname{Re}\,}}

\def\be{\begin{equation}}
\def\ee{\end{equation}}

\newcommand{\ben}{\begin{enumerate}}
\newcommand{\een}{\end{enumerate}}

\newcommand{\blem}{\begin{lem}}
\newcommand{\elem}{\end{lem}}
\newcommand{\bthm}{\begin{thm}}
\newcommand{\ethm}{\end{thm}}
\newcommand{\bcor}{\begin{cor}}
\newcommand{\ecor}{\end{cor}}
\newcommand{\beg}{\begin{exam}}
\newcommand{\eeg}{\end{exam}}
\newcommand{\begs}{\begin{examples}}
\newcommand{\eegs}{\end{examples}}
\newcommand{\bdefe}{\begin{defn}}
\newcommand{\edefe}{\end{defn}}
\newcommand{\bprob}{\begin{prob}}
\newcommand{\eprob}{\end{prob}}
\newcommand{\bques}{\begin{ques}}
\newcommand{\eques}{\end{ques}}
\newcommand{\bei}{\begin{itemize}}
\newcommand{\eei}{\end{itemize}}
\newcommand{\bcon}{\begin{conj}}
\newcommand{\econ}{\end{conj}}
\newcommand{\bop}{\begin{op}}
\newcommand{\eop}{\end{op}}

\newcommand{\bas}{\begin{assertion}}
\newcommand{\eas}{\end{assertion}}

\newcommand{\bfa}{\begin{fact}}
\newcommand{\efa}{\end{fact}}

\newcommand{\bca}{\begin{ca}}
\newcommand{\eca}{\end{ca}}

\newcommand{\bst}{\begin{step}}
\newcommand{\est}{\end{step}}

\newcommand{\bsca}{\begin{sca}}
\newcommand{\esca}{\end{sca}}

\newcommand{\bcl}{\begin{cl}}
\newcommand{\ecl}{\end{cl}}

\newcommand{\bmlem}{\begin{mlem}}
\newcommand{\emlem}{\end{mlem}}

\newcommand{\bscl}{\begin{scl}}
\newcommand{\escl}{\end{scl}}

\newcommand{\bcons}{\begin{conjs}}
\newcommand{\econs}{\end{conjs}}

\newcommand{\bprop}{\begin{prop}}
\newcommand{\eprop}{\end{prop}}

\newcommand{\br}{\begin{rem}}
\newcommand{\er}{\end{rem}}
\newcommand{\brs}{\begin{rems}}
\newcommand{\ers}{\end{rems}}
\newcommand{\bo}{\begin{obser}}
\newcommand{\eo}{\end{obser}}
\newcommand{\bos}{\begin{obsers}}
\newcommand{\eos}{\end{obsers}}
\newcommand{\bpf}{\begin{pf}}
\newcommand{\epf}{\end{pf}}
\newcommand{\ba}{\begin{array}}
\newcommand{\ea}{\end{array}}
\newcommand{\beq}{\begin{eqnarray}}
\newcommand{\beqq}{\begin{eqnarray*}}
\newcommand{\eeq}{\end{eqnarray}}
\newcommand{\eeqq}{\end{eqnarray*}}

\newcommand{\Llra}{\Longleftrightarrow}

\newcommand{\ds}{\displaystyle}

\newcounter{minutes}
\divide\time by 60
\newcounter{hours}
\multiply\time by 60 \addtocounter{minutes}{-\time}

\begin{document}

\bibliographystyle{amsplain}
\title [Generalization of Bohr-type inequality in analytic functions]
{Generalization of Bohr-type inequality in analytic functions}

\def\thefootnote{}
\footnotetext{ \texttt{\tiny File:~\jobname .tex,
          printed: \number\day-\number\month-\number\year,
          \thehours.\ifnum\theminutes<10{0}\fi\theminutes}
} \makeatletter\def\thefootnote{\@arabic\c@footnote}\makeatother

\author{Rou-Yuan Lin}
\address{R.-Y. Lin,  School of Mathematical Sciences, South China Normal University, Guangzhou, Guangdong 510631, China. }
\email{740669790@qq.com}

\author{Ming-Sheng Liu${}^{~\mathbf{*}}$}
\address{M.-S. Liu, School of Mathematical Sciences, South China Normal University, Guangzhou, Guangdong 510631, China.}
\email{liumsh65@163.com}
\author{Saminathan Ponnusamy}\address{S. Ponnusamy, Department of Mathematics,Indian Institute of Technology Madras, Chennai-600 036, India. }\email{samy@iitm.ac.in}

\subjclass[2010]{Primary 30A10, 30C45; Secondary 30H05.}
\keywords{Bohr radius; Bounded analytic function; Alternating series; Bohr inequality. 
 \\
${}^{\mathbf{*}}$ Correspondence should be addressed to Ming-Sheng Liu
}

\begin{abstract}
This paper mainly uses the nonnegative continuous function $\{\zeta_n(r)\}_{n=0}^{\infty}$ to redefine the Bohr radius for
the class of  analytic functions satisfying $\real f(z)<1$   in the unit disk $|z|<1$ and redefine the Bohr radius of the alternating series $A_f(r)$ with analytic functions $f$ of the form $f(z)=\sum_{n=0}^{\infty}a_{pn+m}z^{pn+m}$ in $|z|<1$. In the latter case, one can also get information about Bohr radius for even and odd analytic functions. Moreover, the relationships between the majorant series $M_f(r)$ and the odd and even bits of $f(z)$ are
also established.  We will prove that most of results are sharp.
\end{abstract}

\maketitle
\pagestyle{myheadings}
\markboth{R.-Y. Lin, M-S Liu and S. Ponnusamy}{Generalization of Bohr-type inequality}

\section{Introduction and Preliminaries}\label{LLP-sec1} 

The purpose of the present article is to show several refined and generalized Bohr-type inequality. These results are based on recent investigation and various
generalizations. Moreover, the importance of the classical Bohr inequality is widely appreciated as can be seen from the appearance of many recent investigations on this topic, originated from the geometric function theoretic point of view (cf. \cite{ABS2017,AKP2019,1KP2017,1KP2018}). For a detailed account of
investigation on this topic, one may refer to the articles by Abu-Muhanna et al. \cite{AAP2016}, Defant and Prengel \cite{DePre06},
Garcia et al. \cite{GarMasRoss-2018} and the recent book by Defant et al. \cite{DGMP19}.

\subsection{Classical Bohr inequality and some recent developments}
Let $\mathcal A$ denote the class of all analytic functions defined on the unit disk $\mathbb{D}=\{z\in\mathbb{C}:\,|z|<1\}$.
The present article mainly considers the following two subclasses of $\mathcal A$, namely,
$$\mathcal{B} = \{f\in {\mathcal A}:\, \mbox{ $|f(z)|< 1$ in $\mathbb{D}$}\} ~\mbox{ and }~
\mathcal{P} = \{f\in {\mathcal A}:\, \mbox{ ${\rm Re\,}f(z)< 1$ in $\mathbb{D}$}\}.
$$
While defining the class $\mathcal{B}$, sometimes we consider $|f(z)|\le 1$ (instead of $|f(z)|<1$) which is nothing more than admitting unimodular
constant functions into the set
$\mathcal{B}$. The classical power series inequality for the unit disk due to Harald Bohr in 1914 \cite{B1914} is the following:

\begin{Thm}\label{BohrThA}
If $f\in \mathcal{B}$ and $f(z)=\sum_{n=0}^{\infty}a_nz^n$,  then
\begin{equation*}
M_{f}(r):= \sum_{n=0}^{\infty}|a_n|r^n \le 1 ~\mbox{ for }~|z|=r\leq\frac{1}{3},
\end{equation*}
and the constant $1/3$ is sharp.
\end{Thm}

We observe that $M_{f}(r)$ is rotationally invariant, i.e., $M_{f}(r)=M_{f_\theta}(r)$, where $f_{\theta}(z)=e^{-i\theta}f(z)$ with $\theta=\arg a_0.$
Thus, it suffices to assume in Theorem~A 
that $a_0 \in [0,1).$
	
At first Bohr obtained this inequality only for $r\le 1/6$, but later Wiener, Riesz and Schur independently proved  that  it holds for $r \leq 1/3$ (cf. \cite{PPS2002,S1927,T1962}). 
The number $1/3$ is called the Bohr radius  for $\mathcal{B}$. Other proofs were given in \cite{PS2004,PS2006}. Indeed as pointed out in \cite{PW2020}, in terms of subordination, $f\in \mathcal{B}$ is equivalent to
writing
$$f(z)\prec \varphi_{a_0} (z), \quad  \varphi_ a(z)=\frac{a-z}{1-\overline{ a}z }=a -(1-|a|^2)\frac{z}{1-\overline{a}z },
$$
and thus, from \cite{Rogo-43}, it follows that $|a_n| \leq |\varphi_{a_0} '(0)|=1-|a_0|^2$ for $n\geq 1$  which clearly shows that
$$M_{f}(r) \leq |a_0|+ (1 -|a_0|^2)\frac{r}{1-r}\leq  |a_0|+ (1 -|a_0|)\frac{2r}{1-r} \leq 1   ~\mbox{ for  $r\leq 1/3$}.
$$
The number $1/3$ is optimal because for $a\in (0,1)$, we have
$$M_{\varphi_{a} }(r) = 1 + \frac{1-a}{1-a r}\left ( r(1+2a)-1\right ) >1  ~\Llra~ r >\frac{1}{1+2a},
$$
which by allowing  $a\rightarrow 1^{-}$ shows that the Bohr radius for $\mathcal{B}$ cannot be bigger than $1/3$.

\br\label{LSL1-r1}
We would like to point out that a similar approach as above clearly shows that Theorem~A
continues to hold even 
if we replace $f\in \mathcal{B}$ by $f\in \mathcal{P}$ such that $f(0)\in [0,1)$.
However, an earlier proof of this result may be obtained from \cite[Theorem 2.1]{PPS2002} and  \cite{AAP2016}.
\er

Determination of Bohr radius is being extensively considered nowadays in different settings, 
which includes the discussion on the Bohr radius  for holomorphic functions of several complex variables (cf. \cite{A2000,BK1997,DjaRaman-2000,HHK2009,LL2020}) or  Bloch spaces (cf. Kayumov  et al. \cite[Section~4]{KPS2017}, Liu and Ponnusamy \cite{LLP2018}), Hardy spaces (cf. \cite{BDK5,DjaRaman-2000}) and moreover  the majorant series comparison of two analytic functions satisfying the subordination relation (cf.
\cite{A2010,LPW2020,PVW2021} and the references therein). Kayumov et al. \cite{1KP2017,2KP2018} considered the problem of determining Bohr radius for the classes of odd analytic functions, and for locally univalent functions, and for derivatives of analytic functions,  we refer to \cite{1BD2019,2BD2019}. Recently, it has become very common to apply Bohr radius to harmonic function \cite{EPR2019,1KP2018,KPS2017} and the Bohr radius of harmonic functions is almost parallel to that of the case of  analytic functions.

A version of Bohr's inequality for the family $\mathcal P$ was discussed by Ponnusamy et al. in \cite{PVW2019}. Therefore in Section \ref{LLP-sec3}, we mainly extend the Bohr radius
for the family $\mathcal P$. Moreover, the properties of $M_f(r)$ are extended to the alternating series $A_f(r)$, which was originally investigated in  \cite{ABS2017} by Ali et al. In fact they have shown that if $f\in\mathcal{B}$, then $|A_f(r)|<1$ holds for $r\le 1/\sqrt{3}$. In addition, according to Liu et al. in \cite{LSX2018}, the Bohr radius can be generalized by using the property of $A_f(r)$. Surprisingly, in recent papers, the  authors in \cite{HLP2020,LLP2020,PVW2020,PW2020} refined the Bohr inequality in improved form by replacing $M_{f}(r)$ by $M_{f}(r)+(\frac{1}{1+|a_0|}+\frac{r}{1-r})\Vert f_0\Vert_r^2$ and more recently, the authors in \cite{KKS2021,PVW2021} use nonnegative continuous functions $\{\zeta_n(r)\}_{n=0}^{\infty}$ in $[0,1)$ instead of  $\{r^n\}_{n \geq 0}.$  Refer to Majorant series $M_f(r)$. In Sections \ref{LLP-sec3} and \ref{LLP-sec4} of this paper, refined forms are also used to promote the Bohr phenomenon.


\subsection{Basic notation \label{LSL-subsec1-2}} 
Before we continue the discussion, we fix some notations. Let $f\in \mathcal{A}$,  $f(z)=\sum_{n=0}^{\infty}a_nz^n$  and $f_0(z)=f(z)-f(0)$. Then we define respectively the majorant series and the alternating series of $f$ as follow:
$$M_{f}(r)=|a_0|+M_{f_0}(r) ~\mbox{ and }~  A_{f}(r)=|a_0|+A_{f_0}(r).
$$
where $M_{f_0}(r)$ and $A_{f_0}(r)$ are the majorant series and the alternating series of $f_0$ given by
\begin{equation*}
M_{f_0}(r)=\sum_{n=1}^{\infty}|a_n|r^n  ~\mbox{ and }~  A_{f_0}(r)=\sum_{n=1}^{\infty}(-1)^n|a_n|r^n.
\end{equation*}
We let for convenience,
\begin{equation*}
 \Vert f_0\Vert_r^2=\sum_{n=1}^{\infty}|a_n|^2r^{2n} ~\mbox{ and }~ \Vert f\Vert_r^2 = |a_0|^2+\Vert f_0\Vert_r^2.
\end{equation*}
In addition, if $f\in \mathcal{A}$ and $f(z)=\sum_{n=0}^{\infty}a_{np+m}z^{np+m}$ for some $m\geq 0$ and $p\geq 1$, we define
$$A_{f_{m}}(r)=\sum_{n=1}^{\infty}(-1)^{np+m}a_{np+m}r^{np+m}.
$$
In order to develop the Bohr phenomenon proposed by Kayumov et al. \cite{KKS2021} and developed
further by Ponnusamy et al. \cite{PVW2021}, we need to introduce some new notation. Let $\mathcal{F}$ denote the set of all sequences $\zeta=\{\zeta_n(r)\}_{n=0}^{\infty}$ of nonnegative continuous functions in $[0,1)$ such that the series $\sum_{n=0}^{\infty}\zeta_n(r)$ converges locally uniformly with respect to $r\in[0,1)$. For convenience, we let $\Phi_N(r)=\sum_{n=N}^{\infty}\zeta_n(r)$ whenever $\zeta\in\mathcal{F}$, and in what follows we let
\begin{equation*}
B_N(f, \zeta, r):=\sum_{n=N}^{\infty}|a_n|\zeta_n(r) \ \mbox{for} \ N\ge0;\quad G(f_0, \zeta, r):=\sum_{n=1}^{\infty}|a_n|^2\left(\frac{\zeta_{2n}(r)}{1+|a_0|}+\Phi_{2n+1}(r)\right)
\end{equation*}
for $f(z)=\sum_{n=0}^{\infty}a_nz^n.$ The reason for introducing these quantities are clear from the recent development on this topic (cf. \cite{KKS2021,PVW2021}). In the classical case $\zeta =\{r^n\}_{n\ge 0}$,  $G(f_0, \zeta, r)$ reduces to $\left(\frac{r}{1-r}+\frac{1}{1+|a_0|}\right)\|f_0\|_r^2.$

\subsection{Known important theorems} 
Recently, Kayumov and Ponnusamy et al. \cite{KP2017} (see also \cite{PVW2019}) refined Theorem~A 
in the following form whereas, Liu et al. \cite{LSX2018} proved the estimate of the alternating series.

\begin{Thm}\label{Theo-A} (\cite{KP2017})\, 
Suppose that $f \in \mathcal{B}$ with $f(z)=\sum_{n=0}^{\infty}a_nz^n,$ 
and $S_r$ denotes the area of the image of the subdisk $|z|<r$ under the mapping $f$. Then
\begin{equation*}
M_{f}(r)+\frac{16}{9} \left(\frac{S_r}{\pi} \right)\le 1 \quad for \quad r\le \frac{1}{3},
\end{equation*}
and the constants $1/3$ and $16/9$ cannot be improved. Moreover
\begin{equation*}
M_{f}(r)+|f(z)-a_0|^2 \le 1 \quad for \quad r\le \frac{1}{3},
\end{equation*}
and the constant $1/3$ cannot be improved.
\end{Thm}

In \cite{KKS2021}, a natural generalization of  Theorem~A 
was established. We now recall this result. Further development of this result
for the classes of subordinations and quasi-subordinations may be obtained from \cite{PVW2021}.

\begin{Thm}\label{Theo-B}(\cite{KKS2021})\,
 Let $f\in \mathcal{B}$ with $f(z)=\sum_{n=0}^{\infty}a_nz^n$ and $p\in(0,2]$. If $\zeta_0(r)>\frac{2}{p}\sum_{n=1}^{\infty}\zeta_n(r)$ for $r\in[0,R)$, where $R$ is the minimal positive root of the equation $\zeta_0(x)=\frac{2}{p}\sum_{n=1}^{\infty}\zeta_n(x)$, then the following sharp inequality holds:
\begin{equation*}
B_f(\zeta,p,r):=|a_0|^p\zeta_0(r)+\sum_{n=1}^{\infty}|a_n|\zeta_n(r)\le \zeta_0(r) \ \mbox{ for all $r\le R$.}
\end{equation*}
 In the case when $\zeta_0(x)<\frac{2}{p}\sum_{n=1}^{\infty}\zeta_n(x)$ in some interval $(R,R+\varepsilon)$, the number $R$ cannot be improved.
\end{Thm}

\begin{Thm}\label{Theo-C}(\cite[Theorem 1]{PVW2019})\, 
Let $f\in \mathcal{P}$ with $f(z)=\sum_{n=0}^{\infty}a_nz^n$ and $a_0\in [0,1)$. Then
\begin{equation} \label{TheoC-Eqn}
M_{f}(r)+\left(\frac{1}{1+|a_0|}+\frac{r}{1-r}\right) \|f_{0} \|_{r}^2\le 1
\end{equation}
holds for all $r\le r_*$, where $r_*\approx 0.24683$ is the unique root of the equation $3r^3-5r^2-3r+1=0$ in the interval $[0,1)$.
Moreover, for any  $a_{0}\in (0,1)$ there exists a uniquely defined $r_0=r_0(a_0)\in \left(r_{*},\frac{1}{3}\right)$ such that
$$\sum_{n=0}^\infty |a_n|r^n+\left(\frac{1}{1+a_0}+\frac{r}{1-r}\right)\|f_{0} \|_{r}^2
\le 1
$$
for $r\in [0,r_0]$. The radius $r_0=r_0(a_0)$ can be calculated as the solution of the equation
$$
 \Phi(\lambda,r)= 4r^3\lambda^2- (7r^3+3r^2-3r+1)\lambda +6r^3-2r^2-6r+2 = 0,
$$
where $\lambda = 1 - a_0.$ The result is sharp.
\end{Thm}

At this place, it is worth pointing out that \eqref{TheoC-Eqn} holds for $r \leq 1/3$  (cf. \cite[Theorem 2]{PVW2020}) if the assumption $f\in \mathcal{P}$ is replaced by $f\in \mathcal{B}$.
As mentioned in the introduction, this is in contrast to the case of classical Bohr inequality (refer to Theorem~A 
and Remark \ref{LSL1-r1}).

\begin{Thm}\label{Theo-D}(\cite{LSX2018})\, 
Let $f\in \mathcal{B}$ with $f(z)=\sum_{n=0}^{\infty}a_nz^n$. Then
\begin{equation*}
\left |\ |f(z)|+ A_{f_0}(r)\right | \leq 1~~\mbox{ for  }~~|z|=r \leq\sqrt{2}-1.
\end{equation*}
\end{Thm}

\bprob\label{HLP-prob1} Can we establish an analog of Theorems~B and C
for functions $f\in \mathcal{P}$?
\eprob

The paper is organized as follows. In Section \ref{LLP-sec2}, we mainly state and prove several key lemmas which play key role in the proof of our main Theorems.
In Section \ref{LLP-sec3}, we present partial answers to Problem \ref{HLP-prob1} and present other promotion concerning  functions $f\in \mathcal{P}$.
Moreover, in Section \ref{LLP-sec4}, we also discuss and refine the Bohr-type inequalities of the alternating series $A_f(r)$ and
improve some conclusions of \cite{ABS2017} and \cite{LSX2018}. Throughout the article for $a \in [0,1),$ we let
\begin{equation} \label{automorphism}
\varphi_{a}(z)=\dfrac{a-z}{1-az} ~\mbox{ and }~ \psi _{a} (z)= a-2(1-a)\frac{z}{1-z}, \ z \in \mathbb{D}.
\end{equation}
Note that $\varphi_{a} \in \mathcal{B}$ and $\psi_{a} \in \mathcal{P}$.

\section{Key lemmas and their Proofs}\label{LLP-sec2}

In order to establish our main results, we need the following  lemmas and some important remarks. Two of these lemmas are well-known and the remaining are
new and are key in our discussion.

\begin{Lem}\label{HLP-lem1}{\rm (\cite{AAP2016})}
 If $p\in {\mathcal A}$ and $p(z)=\sum_{n=0}^{\infty}p_nz^n$ such that ${\rm Re}\, p(z)>0$ in $\mathbb{D}$, then we have the sharp inequality
$|p_n| \le 2{\rm Re}\, p_0$ for all $n\ge1$.
\end{Lem}
 Following the notation of Section 1.2, we have

\begin{Lem}\label{HLP-lem3}
{\rm (\cite{PVW2021})}  \quad Let $f\in \mathcal{B}$, $f(z)=\sum_{n=0}^{\infty}a_nz^n$ and $p\in(0,2]$. If $\zeta=\{\zeta_n(r)\}_{n=0}^{\infty}\in\mathcal{F}$, then the following sharp inequality holds:
\begin{equation*}
|a_0|^p\zeta_0(r)+B_1(f,\zeta,r)+G(f_0,\zeta,r)\le|a_0|^p\zeta_0(r)+(1-|a_0|^2)\Phi_1(r)
\end{equation*}
for all $|z|=r\in[0,1)$. The result is sharp for  $\varphi_{a}(z)$ given by \eqref{automorphism}.
\end{Lem} 

\br\label{HLP-re1}
As a  motivation for our discussion, we mention now few useful remarks concerning some special cases of Lemma~G 
for $f\in \mathcal{B}$ with $f(z)=\sum_{n=0}^{\infty}a_nz^n$.
Proof of each item below can be easily derived.

\ben
\item[(i)]  Set $p=1$, $\zeta_{n}(r)=r^{n}$ for $n=km$ for some fixed natural number $k\ge 1$ and $0$, otherwise. Then
Lemma~G  
gives following sharp inequality for $f\in \mathcal{B}$:
\begin{equation}
 \sum_{m=0}^{\infty}|a_{km}|r^{km}+\left (\frac{1}{1+|a_0|}+\frac{r^k}{1-r^k}\right )\|f_0\|_{r}^2 \le |a_0|+(1-|a_0|^2)\frac{r^k}{1-r^k}
\label{liu21}
\end{equation}
for $r\in[0,1)$. Equality is attained for $f(z)= \varphi_{a}(z^k)$, where $\varphi_{a}$ is given by \eqref{automorphism}.
 As
$$|a_0|+(1-|a_0|^2)\frac{r^k}{(1-r^k)} \leq |a_0|+(1-|a_0|)\frac{2 r^k}{1-r^k} \leq 1
$$
for $r \leq \frac{1}{\sqrt[k]{3}}$, the inequality  \eqref{liu21} contains a refined version of \cite[Lemma 2.1]{ABS2017}.
\item[(ii)] For the case $p=1$, $\zeta_{2n-1}(r)=r^{2n-1}~ (n\ge1)$ and $\zeta_{2n}(r)=0 ~(n\ge0)$, we can obtain the following the sharp inequality
\begin{equation}
  \sum_{n=1}^{\infty}|a_{2n-1}|r^{2n-1}\le\frac{r}{1-r^2}(1-\Vert f \Vert_{r}^2) \ \mbox{ for $r\in[0,1)$.}
\label{liu22}
\end{equation}
Equality is attained for $f= \varphi_{a}$, where $\varphi_{a}(z)$ is given by \eqref{automorphism}.

This formula \eqref{liu22} can also be written in another form which establishes the relationship between $M_f(r)$ and $A_f(r)$:
$$M_f(r)\le\frac{2r}{1-r^2}(1-\Vert f\Vert_r^2)+A_f(r) \ \mbox{ for $r\in[0,1)$.}
$$

\item[(iii)] Adding \eqref{liu21} with $k=2$ and \eqref{liu22}, we can easily obtain the following known inequality (which is  \eqref{liu21} with $k=1$):
\begin{equation*}
  M_{f_0}(r)+\left(\frac{1}{1+|a_0|}+\frac{r}{1-r}\right) \|f_0\|_{r}^2 \le (1-|a_0|^2)\frac{r}{1-r} \ \mbox{ for $r\in[0,1)$.}
\end{equation*}
 Equality is attained for $f(z)=\varphi_{a}(z)$, where $\varphi_{a}(z)$ is given by \eqref{automorphism}.
\een
\er

\begin{lem}\label{HLP-lem4}
If $f\in \mathcal{P}$, $a_0=f(0)\in[0,1)$, and $S_r(f)$ denotes the area of the image of the subdisk $|z|<r$ under the mapping $f$, then the following sharp inequality holds
\begin{equation}\label{Lemma4-hypo}
 \frac{S_r(f)}{\pi} \le 4(1-a_0)^2\frac{r^2}{(1-r^2)^2}\quad for \quad |z|=r<1.
\end{equation}
\end{lem}
\bpf
Let $f\in \mathcal{P}$ and $f(z)=\sum_{n=0}^{\infty}a_nz^n$, where $a_0 \in [0,1).$ Applying Lemma~F 
to $p(z)=1-f(z)$, we easily have
$|a_n|\le2(1-a_0)$ for $n\ge 1$ and thus,
$$\frac{S_r(f)}{\pi} =\frac{1}{\pi}\iint_{\mathbb{D}}|f'(z)|^2\,dxdy =\sum_{n=1}^{\infty}n|a_n|^2r^{2n}  \le4(1-a_0)^2\frac{r^2}{(1-r^2)^2}
$$
 and the inequality  \eqref{Lemma4-hypo} follows. It is easy to see that the inequality \eqref{Lemma4-hypo} becomes an equality for the
function $f(z)=\psi_{a}(z)$, where $\psi_{a}(z)$ is given by \eqref{automorphism}. This is because for this function, $a_0=a$ and $a_n=-2(1-a)$  for $n \geq 1$.
\epf

\begin{lem}\label{HLP-lem5} Suppose that  $f \in \mathcal{B}$ with  $f(z)=\sum_{n=0}^{\infty}a_{np+m}z^{np+m}$, where $p\in \mathbb{N}$ and $0\le m < p$. Then
\begin{equation}
 |f(z)|+ r^m A(r)  \le 1
\label{liu23}
\end{equation}
for $|z|=r\le r_{p,m}$, where $r_{p,m}$ is the minimal positive root of the equation
$$r^{3p-m}-2r^{3p}-3r^{2p}-r^{p-m} + 1=0,
$$
and
$$A(r) = \sum_{n=1}^{\infty}|a_{2np+m}|r^{2np}+ \left(\frac{1}{1+|a_m|}+\frac{r^{2p}}{1-r^{2p}} \right)\sum_{n=1}^{\infty}|a_{np+m}|^2r^{2np}.
$$
The radius $r_{p,m}$ is best possible.
\end{lem}

\bpf 
We write $f$ as $f(z)=z^mg(z^p)$, where $p\geq 1$ and $m >0$ are fixed and $g(z)=\sum_{n=0}^{\infty}b_nz^{n}\in\mathcal{B}$ with $b_n=a_{np+m}\ (n \geq 0)$. Schwarz-Pick Lemma
applied to $g$ gives
\begin{equation*}
|f(z)|\le r^m\left( \frac{r^p+|a_{m}|}{1+r^p|a_{m}|} \right) \ \mbox{for all $z=re^{i\theta}\in \ID$}.
\end{equation*}
This inequality and Lemma \ref{HLP-th5}(i) (with $\zeta_{2n}(r)=r^{2n}$, $n\ge0$) proved later give, for $r\in[0,1)$,
\begin{eqnarray}
  |f(z)|+ r^m A(r)  &\le& r^m\left[\frac{r^p+|a_{m}|}{1+r^p|a_{m}|}  +(1-|a_m|^2)\frac{r^{2p}}{1-r^{2p}} \right ] \label{liu25}\\
  & = &1+\dfrac{P_{r}(|a_m|)}{(1-r^{2p})(1+r^p|a_m|)} \nonumber
\end{eqnarray}
where
\begin{eqnarray*}
  &&P_r(x)=(xr^{3p}-x^3r^{3p+m}+xr^{3p+m}-r^{3p+m})\nonumber\\
  &&\hspace{1.5cm}+(r^{2p+m}+r^{2p}-xr^{2p+m}-x^2r^{2p+m})+(r^{p+m}-xr^p)+(xr^m-1).
\end{eqnarray*}
 Thus, to prove that the left hand side of (\ref{liu25}) is smaller than or equal to $1$, it suffices to show that $P_r(x)\leq 0$ for $x\in [0,1)$.  Note that
$$P_r(1)=(1-r^m)(r^{3p}+r^{2p}-r^{p}-1)
 =(1-r^m)(1+r^p)(r^{2p}-1)\leq  0.
$$
We claim that $P_r(x)$ is an increasing function of $x$ for $x\in [0,1]$. 
Now a direct computation shows that
\begin{equation*}
P_r'(x)=r^{3p}-3x^2r^{3p+m}+r^{3p+m}-r^{2p+m}-2 xr^{2p+m}-r^p+r^m,
\end{equation*}
and
$P_r''(x)=-6xr^{3p+m}-2r^{2p+m}\le0,
$
for all $x\in[0,1]$.  Thus, $P_r'$ is a decreasing function of $x \in [0,1].$  Consequently, by the condition of the theorem,  we have
\begin{equation*}
P_r'(x)\ge P_r'(1)=r^m[r^{3p-m}-2r^{3p}-3r^{2p}-r^{p-m}+1]\ge0 ~\mbox{ for $r\le r_{p,m}$.}
\end{equation*}
Thus for each $r\le r_{p,m}$, $P_r(x)$ is an increasing function of $x$ on $[0,1]$, which implies that $P_r(x)\le P_r(1)<0$ for all
$x\in[0,1]$  and $r\le r_{p,m}$.

To show that the radius $r_{p,m}$ is best possible, we let $a\in[0,1)$ and consider the function
\begin{eqnarray*}
f(z)=z^m \left(\frac{a+z^p}{1+az^p} \right)=az^m+(1-a^{2})\sum_{n=1}^{\infty}(-1)^{n-1}a^{n-1}z^{np+m},\quad  z\in \ID.
\end{eqnarray*}
A simple computation gives
$$ A(r) = (1-a^2)\frac{ar^{2p}}{1-a^2r^{2p}} +\left(\frac{1}{1+a}+\frac{r^{2p}}{1-r^{2p}} \right)(1-a^2)^2\frac{r^{2p}}{1-a^2r^{2p}} =(1-a^{2})\frac{r^{2p}}{1-r^{2p}}
$$
and thus,  we find that
\begin{eqnarray*}
|f(r)|+r^mA(r)
&=&r^m\left[\frac{a+r^p}{1+ar^p}  +(1-a^{2})\frac{r^{2p}}{1-r^{2p}} \right].
\end{eqnarray*}
Comparison of this expression with the right hand side of the inequality in  (\ref{liu25}) delivers the asserted sharpness.
\epf

\begin{lem} \label{HLP-lem6} For $p>0$, let
$$T(x)=\frac{1-x^p}{1-x} ~\mbox{ for $x\in [0,1)$.}
$$
Then
\begin{enumerate}
\item[{\rm (i)}] for each $p\in(0,1],$  we have	$T(x)\ge p$ for all $x\in [0,1)$.

\item[{\rm (ii)}] for each $p>1,$  we have $1=T(0)\leq T(x) \leq \lim_{x \to 1^{-}}T(x)=p$.
\end{enumerate}
	
\end{lem}
\bpf
For $y \in [0,1]$ and $0<p\leq1,$ the Binomial theorem quickly gives $(1-y)^p \leq 1-py$ so that (i) follows if we replace $y$ by $1-x.$ This fact and the second part also follow from the fact that
$f(y)=(1-y)^p-(1-py)
$
is decreasing function of $y,$ when $0<p\leq 1$ and is an increasing function of $y,$ when $p>1.$
\epf

\section{Bohr type inequalities of analytic functions $f$ with $\real f(z) < 1$}\label{LLP-sec3}

We first present Bohr-type radius for $f \in \mathcal{P}$ instead of $f \in \mathcal{B}$ and cover a more general  Bohr-type sum.
\bthm\label{HLP-th1}
 Suppose that $f \in \mathcal{P}$ and $f(z)=\sum_{n=0}^{\infty}a_nz^n$ with $a_0\in[0,1)$. 
Then we have the following:
\begin{itemize}
\item[(i)] If $p\in(0,1]$ and $\{\zeta_n(r)\}_{n=0}^{\infty}\in\mathcal{F}$ satisfies
\beq \label{thm1-eqn1}
\zeta_0(r)\ge\frac{2}{p} \Phi_1(r)
\eeq
for $r\in[0,R_1]$, where $\Phi_1(r)=\sum_{n=1}^{\infty}\zeta_n(r)$ and $R_1$ is the minimal positive root of the equation
\begin{equation}
\zeta_0(r)=\frac{2}{p}\Phi_1(r)
\label{liu31}
\end{equation}
then the following sharp inequality holds:
\begin{equation}
B_f(\zeta,p,r):= a_0^p\zeta_0(r)+\sum_{n=1}^{\infty}|a_n|\zeta_n(r)\le\zeta_0(r)\ \mbox{for all} \ r\le R_1.
\label{liu32}
\end{equation}
In the case when $\zeta_0(r)<\frac{2}{p}\Phi_1(r)$ 
 in some interval $(R_1,R_1+\varepsilon)$, the number $R_1$ cannot be improved.

\item[(ii)] If $p>1$, $p\in\mathbb{N}$ and
$$
\zeta_0(r)\ge\frac{2}{1+a_0+\cdots+a_0^{p-1}}\Phi_1(r)
$$
for $r\in[0,R_p]$, where $R_p$ is the minimal positive root of the equation $$\zeta_0(r)=\frac{2}{1+a_0+\cdots+a_0^{p-1}}\Phi_1(r),
$$ then we have the sharp inequality
\begin{equation*}
B_f(\zeta,p,r)\le\zeta_0(r) \quad for \quad all \quad r\le R_p.
\end{equation*}
In particular, when $p=2$, $B_f(\zeta,2,r)\le\zeta_0(r)$ for all $r\le R_2$, where $R_2$ is the root of $\zeta_0(r)=\frac{2}{1+a_0}\Phi_1(r)$
\end{itemize}
\ethm

\bpf
(i) As $f \in \mathcal{P},$ we have $\real (1-f(z))>0$ in $\mathbb{D}$ and so by Lemma~F, 
we have $|a_n|\le 2(1-a_0)$ for all $n\ge1$ and thus, we see that
\beq
B_f(\zeta,p,r)&\le&a_0^p\zeta_0(r)+2(1-a_0)\sum_{n=1}^{\infty}\zeta_n(r)=a_0^p\zeta_0(r)+2(1-a_0)\Phi_1(r)\nonumber\\
&=&\zeta_0(r)+2(1-a_0)\Bigg[\Phi_1(r)-\frac{1-a_0^p}{2(1-a_0)} \zeta_0(r)\Bigg]\label{LSL-eq1b}\\
&\le&\zeta_0(r)+2(1-a_0)\Bigg[\Phi_1(r)-\frac{p}{2}\zeta_0(r)\Bigg]\nonumber\\
&\le&\zeta_0(r),\ \nonumber \mbox{by Lemma \ref{HLP-lem6}(i) and \eqref{thm1-eqn1}, }
\eeq
for all $r\le R_1$. This proves the inequality (\ref{liu32}).

Now, let us prove that $R_1$ is an optimal number. For the function $f(z)=\psi_a (z)$, where $\psi _a$ is defined by
\eqref{automorphism}, we easily see that
\be\label{LSL-eq1a}
B_f(\zeta,p,r)
=\zeta_0(r)+2(1-a)\Bigg[\Phi_1(r)-\left(\frac{1-a^p}{2(1-a)}\right)\zeta_0(r)\Bigg],
\ee
which is bigger than $\zeta_0(r)$ if
$$
\Phi_1(r)-\left(\frac{1-a^p}{2(1-a)}\right)\zeta_0(r)>0.
$$
In fact, for $r>R_1$ and $a$ close to $1$, we see that
\begin{equation*}
\lim_{a\rightarrow1^-}\bigg[\Phi_1(r)
-\left(\frac{1-a^p}{2(1-a)}\right)\zeta_0(r)\bigg]=\Phi_1(r)
-\frac{p}{2}\zeta_0(r)\geq0,
\end{equation*}
showing that the number $R_1$ is  best possible.

(ii) When $p\in\mathbb{N}\setminus\{1\}$, we easily have by \eqref{LSL-eq1b} that
\begin{eqnarray*}
B_f(\zeta,p,r)&\le&\zeta_0(r)+2(1-a_0)\bigg[\Phi_1(r)
-\frac{1+a_0+\cdots+a_0^{p-1}}{2}\zeta_0(r)\bigg]
\end{eqnarray*}
so that $B_f(\zeta,p,r)\le\zeta_0(r)$, only when $\zeta_0(r)\ge\frac{2}{1+a_0+\cdots+a_0^{p-1}}\Phi_1(r)$.

To prove that $R_p$ is an optimal number, we consider the function $f(z)=\psi_a (z)$, where $\psi _a$ is defined by
\eqref{automorphism}. For this function, by \eqref{LSL-eq1a} we get that
$$B_f(\zeta,p,r)=\zeta_0(r)+2(1-a)\bigg[\Phi_1(r)
-\left(\frac{1+a+\cdots+a^{p-1}}{2}\right)\zeta_0(r)\bigg],
$$
and the sharpness follows as with the previous arguments.
\epf

\br\label{HLP-re2}
\begin{enumerate}
\item[(1)] For  $p>1$, according to the proof of Lemma \ref{HLP-lem6}(ii), 
Theorem \ref{HLP-th1}(i) holds by replacing the factor $\frac{2}{p}$ in (\ref{liu31}) by $2$ and also at other three places.
 However, the result will not be good for $p\ge2$.

\item[(2)] Clearly, Theorem \ref{HLP-th1}(i) for $p=1$ and $\zeta_n(r)=r^n ~(n\ge0)$, gives the classical  Bohr radius $1/3$  for $\mathcal{P}$. See also
Remark \ref{LSL1-r1}, \cite{AAP2016} and \cite[Theorem 2.1]{PPS2002}.
\begin{itemize}
\item[(i)] If $f \in \mathcal{P}$ with $a_0=f(0)\in [0,1)$, $p=2$ and $\zeta_n(r)=r^n ~(n\ge0)$, then
$$a_0^2+\sum_{n=1}^{\infty}|a_n|r^n\le1 \ \mbox{for}  \ r\le R_2=\frac{1+a_0}{3+a_0}\in [1/3,1/2).
$$
 Note that in the case of $f \in \mathcal{B},$ it is well-known that the corresponding sharp value of $R_2$ is $1/2.$
\item[(ii)] If $f \in \mathcal{P}$ with $a_0=f(0)\in [0,1)$,  $p\in(0,1]$, and
$$\zeta_{n}(r)=\left\{\begin{array}{ll}
r^n & \mbox{ if $n=km,\ m \geq1,$}\\
0 & \mbox{ otherwise,}
\end{array}\right.
$$ where $k \in \mathbb{N}$ is fixed, then we get the following sharp inequality:
\begin{equation*}
a_0^p+\sum_{m=1}^{\infty}|a_{km}|r^{km}\le1 \ \mbox{for} \ r\le \sqrt[k]{\frac{p}{2+p}},\ k \geq 1.
\end{equation*}
This result is same as for $f\in\mathcal{B}$. See also \cite[Remark 1]{PVW2019} for a refined version (when $k=1$).
\end{itemize}
\end{enumerate}
\er

\bcor\label{HLP-cor1}
  Suppose that $f \in \mathcal{P}$ and $f(z)=\sum_{n=0}^{\infty}a_nz^n$ with  $a_0\in[0,1)$. 
Then the following sharp inequality holds:
\begin{equation*}
B_f(\zeta,1,r)+G(f_0,\zeta,r):= a_0\zeta_0(r)+\sum_{n=1}^{\infty}|a_n|\zeta_n(r)+\sum_{n=1}^{\infty}|a_n|^2\left(\frac{\zeta_{2n}(r)}{1+|a_0|}+\Phi_{2n+1}(r)\right)\le\zeta_0(r)
\end{equation*}
for all $r\in[0,r_1)$, where $r_1$ is the minimal positive root of the equation
\begin{equation*}
 \zeta_0(r)=2\Phi_1(r)+4\sum_{n=1}^{\infty}\Phi_{2n}(r).
\end{equation*}
  Here $G(f_0,\zeta,r)$ is defined in Section \ref{LSL-subsec1-2} and $B_f(\zeta,1,r)$ is given by \eqref{liu32}.
\ecor

\bpf
The proof 
follows by using the method of proof of Theorem \ref{HLP-th1}(i). 
 Note that when $\zeta_n(r)=r^n$ $(n\ge0)$, Corollary \ref{HLP-cor1} is same as the first part of  Theorem~D. 
\epf

The following results establish analogs of Theorems~B and C
when $f \in \mathcal{P}.$  

\bthm\label{HLP-th2}
Suppose that $f \in \mathcal{P}$ and $f(z)=\sum_{n=0}^{\infty}a_nz^n$ with  $a_0\in[0,1)$.  Also, let $q\ge1$, $m\in\mathbb{N} $,
$\{\zeta_n(r)\}_{n=0}^{\infty}\in\mathcal{F}$ and
$$E_f(\zeta,p,r):=a_0^p\zeta_0(r)+\sum_{n=1}^{\infty}|a_n|\zeta_n(r)+|f(z^m)-a_0|^q.
$$
Then we have the following:
\begin{enumerate}
\item[{\rm (i)}] If $p\in(0,1]$ and
\beq\label{thm2-**}
\zeta_0(r)\ge\frac{2}{p}\bigg[\Phi_1(r)+2^{q-1} \frac{r^{mq}}{(1-r^m)^q} \bigg]
\eeq
for $r\in[0,R_{m,q}^p]$, where $\Phi_1(r)=\sum_{n=1}^{\infty}\zeta_n(r)$ and  $R_{m,q}^p$ is the minimal positive root of the equation
\begin{equation}
\zeta_0(r)=\frac{2}{p}\bigg[\Phi_1(r)+2^{q-1} \frac{r^{mq}}{(1-r^m)^q} \bigg],
\label{liu33}
\end{equation}
then
\begin{equation}
E_f(\zeta,p,r)\le\zeta_0(r)\ \mbox{for all}  \ r\le R_{m,q}^p.
\label{liu34}
\end{equation}

When $p=1$, the result is sharp for all $q\ge1$ and $m\in\mathbb{N}$. Moreover, when $q=1$, this result is sharp for all $p\in(0,1]$ and $m\in\mathbb{N}$.

\item[{\rm (ii)}] If $p=2$ and $\zeta_0(r)\ge\frac{2\Phi_1(r)+2^q(1-a_0)^{q-1}\frac{r^{mq}}{(1-r^m)^q}}{1+a_0}=:\psi(r)$, then
\begin{equation*}
E_f(\zeta,2,r)\le\zeta_0(r) \ \mbox{for all} \ r\le R_{m,q}^2.
\end{equation*}
where $R_{m,q}^2$ is the minimal positive root of the equation $\zeta_0(r)=\psi(r)$, and the result is sharp.
\end{enumerate}
\ethm
\bpf
(i) Let $0<p\leq1.$ Then the triangle inequality and assumption that $f \in \mathcal{P}$ gives
\begin{equation*}
|f(z^m)-a_0|^q=\bigg|\sum_{n=1}^{\infty}a_nz^{mn}\bigg|^q\le \left(\sum_{n=1}^{\infty}|a_n|r^{mn}\right)^q \leq 2^q(1-a_0)^q\left(\frac{r^{m}}{1-r^m}\right)^q
\end{equation*}
and thus, Lemma~F 
yields that
\beq\label{thm2***}
E_f(\zeta,p,r)&\le& a_0^p\zeta_0(r)+2(1-a_0)\Phi_1(r)
+2^q(1-a_0)^q \frac{r^{mq}}{(1-r^m)^q} \\
&=& \zeta_0(r)+2(1-a_0)h(p,a_0,r), \nonumber
\eeq
where
$$
h(p,a_0,r)=\Phi_1(r)+2^{q-1}(1-a_0)^{q-1} \frac{r^{mq}}{(1-r^m)^q} - \frac{1-a_0^p}{2(1-a_0)} \zeta_0(r),
$$
and, by Lemma \ref{HLP-lem6}, it follows that
$$
h(p,a_0,r) \leq H(p,a_0,r):=\Phi_1(r)+2^{q-1}(1-a_0)^{q-1} \frac{r^{mq}}{(1-r^m)^q} -\frac{p}{2}\zeta_0(r).
$$
We need to prove that $E_f(\zeta,p,r) \leq \zeta_0(r)$ and for this it suffices to show that $H(p,a_0,r) \leq 0.$

As $a_0 \in [0,1)$ and $q\ge1$, it is clear that $H(p,a_0,r)$ is a decreasing function of $a_0$
and thus,
\begin{equation*}
H(p,a_0,r)\le H(p,0,r)=\Phi_1(r)+2^{q-1} \frac{r^{mq}}{(1-r^m)^q} -\frac{p}{2}\zeta_0(r)
\end{equation*}
which is non-positive whenever \eqref{thm2-**} holds.
Thus, (\ref{liu34}) holds.

Now let us prove that the radius $R_{m,q}^p$ is sharp for both cases. For $a \in [0,1),$ consider the function
$f(z)=\psi_a (z)$, where $\psi _a$ is defined by \eqref{automorphism}. For this function, it is a simple exercise to see that
\begin{equation}\nonumber
\begin{split}
E_f(\zeta,p,r)&= a^p\zeta_0(r)+2(1-a)\Phi_1(r)
+2^q(1-a)^q \frac{r^{mq}}{(1-r^m)^q} \\
&= \zeta_0(r)+2(1-a)G(a,p,r)
\end{split}
\end{equation}
which is bigger than $\zeta_0(r)$ if
$$G(a,p,r):=\Phi_1(r)+2^{q-1}(1-a)^{q-1} \frac{r^{mq}}{(1-r^m)^q} - \frac{1-a^p}{2(1-a)} \zeta_0(r)>0.
$$
\textbf{Case 1:} $p=1.$ For $a\rightarrow0 ^{+}$ and $r>R_{m,q}^1$, we see that
$$
\lim_{a\rightarrow 0 ^{+}}G(1,a,r)
= \Phi_1(r)+2^{q-1} \frac{r^{mq}}{(1-r^m)^q} -\frac{1}{2}\zeta_0(r)
>0,
$$
 (since $\zeta_0(r)<2[\Phi_1(r)+2^{q-1} \frac{r^{mq}}{(1-r^m)^q}]$ is equivalent to $r>R_{m,q}^1$)  showing that the number $R_{m,q}^1$ is best possible.\\

\noindent
\textbf{Case 2:} $q=1$. For $a\rightarrow 1^{-}$ and $r>R_{m,1}^p$, we see that
$$\lim_{a\rightarrow1^{-}}G(p,a,r)\geq \Phi_1(r)+\frac{r^m}{1-r^m}-\frac{p}{2}\zeta_0(r)>0
$$
 (since $\zeta_0(r)<\frac{2}{p}[\Phi_1(r)+\frac{r^{m}}{(1-r^m)}]$ is equivalent to $r>R_{m,1}^p$)  showing that the number $R_{m,1}^p$ is best possible.
\vskip 1.5mm

(ii)  $p=2$ and $q\ge1.$\\
 In this case, \eqref{thm2***} gives
$
E_f(\zeta,2,r)\le\zeta_0(r)+2(1-a_0)I(2,a,r),
$
where
$$I(2,a,r)=\Phi_1(r)+2^{q-1}(1-a_0)^{q-1}\frac{r^{mq}}{(1-r^m)^q}-\frac{1+a_0}{2}\zeta_0(r).
$$
Thus, $E_f(\zeta,2,r)\le\zeta_0(r)$, whenever $I(2,a,r)\le0$ which is equivalent to the condition
stated in the statement of the theorem, namely, the condition $\zeta_0(r) \ge \psi (r)$.

Finally, to prove that $R_{m,q}^2$ is an optimal number, we just consider the function $\psi _a$ defined by \eqref{automorphism}, and follow the usual procedure.
This completes the proof.
\epf

\br\label{HLP-re3}
\begin{enumerate}
\item[{\rm (1)}]
Obviously in the limiting case $m\rightarrow{\infty}$, Theorem \ref{HLP-th2} reduces to Theorem \ref{HLP-th1}(i). Moreover, when $p\in(0,1]$,
an appropriate choices for $\zeta_{n}(r)$, the limiting case
$m\rightarrow{\infty}$ of Theorem \ref{HLP-th2} also contains Remark \ref{HLP-re2}(ii).

\item[{\rm (2)}]
When $p=1$ and $\zeta_n(r)=r^n \ (n\ge0)$, only when $m$ is sufficiently large, we have $r\le R^1_{m,q}$ with $R^1_{m,q} \approx 1/3.$

\item[{\rm (3)}]
 If $\zeta_n(r)=r^n$, then, from   Table \ref{Tab1}, we see that $R_{m,q}^p$ gradually tends to $1/3$ in the case of $p=1$, and
tends to $1/5$ in the case $p=1/2$, when $m$ increases. 
\begin{table}[htbp]
\centering
\caption{Relevant results obtained by Theorem \ref{HLP-th2}}
\begin{tabular}{cccccccccccc}
\hline
$m$ &$R_{m,2}^1$ & $m$ & $R_{m,2}^1$ & $m$ & $R_{m,1}^1$ &$m$ & $R_{m,1}^1$ & $m$ & $R_{m,1}^{\frac{1}{2}}$ &$m$ & $R_{m,1}^{\frac{1}{2}}$\\
\hline
1 &0.236068& 2& 0.321336& 1& 0.200000& 2&0.289898 &1 &0.111111& 2& 0.178395 \\
\hline
3 &0.332047& 4& 0.333195& 3& 0.318201& 4&0.328083 & 3 &0.195177& 5& 0.199796 \\
\hline
5 &0.333318& 7& 0.333333& 5& 0.331541& 10&0.333326 &10 &0.199999 & 30& 0.199999 \\
\hline
10 &0.333333& 15& 0.333333& 15& 0.333333& 20&0.333333 &50 &0.2& 60& 0.2\\
\hline
\end{tabular}\label{Tab1}
\end{table}
\end{enumerate}
\er

\bcor\label{HLP-cor2}
  Suppose that $f \in \mathcal{P}$ and $f(z)=\sum_{n=0}^{\infty}a_nz^n$ with  $a_0\in[0,1)$. Then
\begin{equation*}
a_0+\sum_{n=1}^{\infty}|a_n|r^n+|f(z)-a_0|^2 \le 1 \ \mbox{ for  $|z|=r\le\sqrt{5}-2$,}
\end{equation*}
and the result is sharp for the function $\psi _a$ defined by \eqref{automorphism}. Moreover,
\begin{equation*}
a_0^2+\sum_{n=1}^{\infty}|a_n|r^n+|f(z)-a_0|^2 \le 1 \ \mbox{for $|z|=r\le\frac{1}{3}$}
\end{equation*}
if and only if $a_0\in[\frac{1}{2},1).$ The result is sharp for the function $\psi _a$ defined by \eqref{automorphism}.
\ecor

The following result is a generalization of Theorems~B and C 
for $\mathcal{P}$.

\bthm\label{HLP-th3}
 Suppose that $f \in \mathcal{P}$, $\lambda > 0$ and $f(z)=\sum_{n=0}^{\infty}a_nz^n$ with  $a_0\in[0,1)$. If $\{\zeta_n(r)\}_{n=0}^{\infty}\in\mathcal{F}$ satisfies the inequality
$$\zeta_0(r)\ge\frac{2}{1+a_0} \left [\Phi_1(r)+2\lambda (1-a_0)\frac{r^2}{(1-r^2)^2}\right ] =:\phi(a_0,r),
$$
then the following sharp inequality holds:
$$a_0^2\zeta_0(r)+\sum_{n=1}^{\infty}|a_n|\zeta_n(r)+\lambda\left(\frac{S_r(f)}{\pi}\right)\le\zeta_0(r)\ \mbox{for all} \ r\le R_{\lambda,2},
$$
where $R_{\lambda,2}$ is the minimal positive root of the equation $\zeta_0(r)=\phi(a_0,r)$. This result is sharp for the function $\psi _a(z)$ defined by
\eqref{automorphism}.
\ethm

\bpf
 With the help of Lemma \ref{HLP-lem4}, the proof is similar as in Theorem \ref{HLP-th2}(ii),  and so it will not be repeated. Moreover, to prove that $R_{\lambda,2}$ is an optimal number, we consider the function $\psi _a$ defined by \eqref{automorphism},
\epf

For $f \in \mathcal{P},$ because $a_0^2+\sum_{n=1}^{\infty}|a_n|r^n\le1 \ \mbox{holds for all} \ r\le 1/3$, 
it is natural to state the following refined form. In fact as a consequence of Theorem \ref{HLP-th3}, we have the following. 

\bcor\label{HLP-cor3}
Suppose that $f \in \mathcal{P},$ and $f(z)=\sum_{n=0}^{\infty}a_nz^n$ with $a_0\in(0,1)$. Then
\begin{equation*}
a_0^2+\sum_{n=1}^{\infty}|a_n|r^n+\frac{16 a_0}{9(1- a_0)}\left(\frac{S_r}{\pi}\right) \le 1 \quad \mbox{for}\quad r\le \frac{1}{3},
\end{equation*}
and the constants $1/3$ and $\frac{16 a_0}{9(1- a_0)}$ cannot be improved, as shown by the function $\psi _a$ defined by \eqref{automorphism}.
\ecor

Finally, we present an improved version of Theorem~D. 

\bthm\label{HLP-th4}
Suppose that $f \in \mathcal{P},$ $f(z)=\sum_{n=0}^{\infty}a_nz^n$ with $a_0\in[0,1)$, and let
$$D_{\lambda}(r)=a_0+\sum_{n=1}^{\infty}|a_n|r^n+\left(\frac{1}{1+a_0}+\frac{r}{1-r}\right)\sum_{n=1}^{\infty}|a_n|^2r^{2n}+\lambda \left(\frac{S_r}{\pi}\right).
$$
Then
\begin{equation}
D_{\lambda}(r) \le 1 \ \mbox{for $\ds |z|=r\le \frac{1}{5-2a_0}$},
\label{liu37}
\end{equation}
 where $\lambda=  8/9\approx0.888889 .$ The constants  $1/(5-2a_0)$ and $\lambda=8/9$ cannot be improved.
\ethm

\bpf
Clearly, $D_{\lambda}(r)$ is an increasing function of $r$ and hence, we only need to prove the inequality \eqref{liu37} for $r=\rho,$
where   $\rho=1/(5-2a_0)$. As before, by Lemma \ref{HLP-lem4} and the fact that $|a_n|\leq 2(1-a_0)$ for $n \geq 1,$ we deduce that
\begin{equation}\nonumber
\begin{split}
D_{\lambda}(\rho)&=a_0+\sum_{n=1}^{\infty}|a_n|\rho^n+\left(\frac{1}{1+a_0}+\frac{\rho}{1-\rho}\right)\sum_{n=1}^{\infty}|a_n|^2\rho^{2n}+ \lambda \left (\frac{S_{\rho}}{\pi} \right )\\
&\le a_0+  \frac{1-a_0}{2-a_0}+\left(\frac{1}{1+a_0}+\frac{1}{4-2a_0}\right)\frac{(1-a_0)^2}{(2-a_0)(3-a_0)} \\
&\quad +\lambda \frac{(1-a_0)^2(5-2a_0)^2}{4(2-a_0)^2(3-a_0)^2}\\
& =1+\Phi_4(a_0),
\end{split}
\end{equation}
where
\beqq
\Phi_4(x)&=&(1-x)\bigg[\frac{1}{2-x}+\bigg(\frac{1}{1+x}+\frac{1}{4-2x}\bigg)\frac{1-x}{(2-x)(3-x)}+\lambda \frac{(1-x)(5-2x)^2}{4(2-x)^2(3-x)^2}-1\bigg]\\
&& =\frac{1-x}{4(2-x)^2(3-x)^2(1+x)}L_\lambda (x).
\eeqq
Here (we will replace $\lambda$ with $\mu$ in order to determine the value of $\lambda$) and below:
\begin{equation*}
L_\mu (x)=-4x^5+32x^4-82x^3+58x^2+38x-42-4\mu x^4+20\mu x^3-21\mu x^2-20\mu x+25\mu.
\end{equation*}

Now  $D_{\lambda}(\rho)\le1$ provided $L_\lambda (x)\le0$ (that is $L_\mu (x)\le0$) for all $x\in[0,1]$. Obviously,
$$
L_\mu(1)=0, ~\mbox{ and }~ L_\mu(0)=  -42+25\mu \le 0 ~\mbox{ for $\mu\le 1.68$} .
$$
Next, we find that
\begin{equation*}
L'_{\mu}(x)=-20x^4+128x^3-246x^2+116x+38-16\mu x^3+60\mu x^2-42\mu x-20\mu,
\end{equation*}
and therefore, $L'(1)=16-18\mu$ showing that
$L_\mu'(1)\ge0$ for $\mu \le \frac{8}{9} $ and $L_\mu'(1)\le0$ for $\mu\ge 8/9$.

Because $L_{\mu}(1)=0$ and $L_{\mu}(x)$ is a differentiable function, so that if $L'_{\mu}(1)<0$, then there exists an $x_0\in (1-\varepsilon, 1)$ such that $L_{\mu}(x_0)>0$,
which does not meet the requirements. Therefore, we should satisfy $L'_{\mu}(1)\geq 0$, that is, $\mu\le 8/9 $.

When $\mu\in[0,8/9]$, we have $L_{\mu}(0)=-42+25\mu\leq 0$ and $L_{\mu}(1)=0$. Since
\begin{equation*}
\frac{\partial L'_{\mu}(x)}{\partial\mu}=-16x^3+60x^2-42x-20\leq 0\quad \mbox{ for } x\in [0, 1],
\end{equation*}
so that for $\mu\in[0,\frac{8}{9}]$, we have
\begin{equation*}
L'_{\mu} (x)\ge-20x^4+128x^3-\frac{128}{9}x^3-246x^2+\frac{160}{3}x^2+116x-\frac{112}{3}x+38-\frac{160}{9}:=l(x)
\end{equation*}
With the aid of Mathematica,
we know that $l(0)\approx20.222222$, $l(1)=0$, $l(x)$ has no roots in the interval $[0,1)$ and $l(x)$ is continuous.
Then $L'_{\mu}(x)\ge l(x)>0$ for all $x\in[0,1)$. This implies that $L_{\mu}(x)\leq L_{\mu}(1)=0$ in $[0,1]$, that is, $D_{\mu}(\rho)\leq 1$ for $\mu\in[0,8/9]$.

In conclusion, the maximum value that $\mu$ can taken as $8/9$. Thus $L_{\lambda}(x)\le0$ for $\lambda =8/9$ and hence for all $0<\lambda\le 8/9.$

To prove that the constant $\lambda=8/9$ is sharp, we consider the function $\psi _a$ defined by \eqref{automorphism}.
For this function, we compute
\begin{equation*}
D_{\lambda_{1}}(\rho )=a_0+\sum_{n=1}^{\infty}|a_n|\rho ^n+\left(\frac{1}{1+a_0}+\frac{\rho}{1-\rho}\right)\sum_{n=1}^{\infty}|a_n|^{2}\rho^{2n}+\lambda_1\left(\frac{S_\rho}{\pi}\right).
\end{equation*}
which simplifies to
\begin{equation*}
D_{\lambda_{1}}(\rho )=1+(1-a)\frac{L_{\lambda_{1}}(a)}{4(2-a)^2(3-a)^2(1+a)},
\end{equation*}
where $L_{\lambda_{1}}(a)$ is a differentiable function and defined as above.
For $\lambda_{1}>\lambda=8/9,$ we get that $L'_{\lambda_{1}}(1)<0$  and by Mathematica has $L_{\lambda_{1}}(0.999999)>0$, which implies $D_{\lambda_{1}}(\rho)>1$.
This proves the sharpness assertion and the proof of Theorem \ref{HLP-th4} is complete.
\epf

\br\label{HLP-re5}
Theorem \ref{HLP-th4} can be regarded as an extension of Theorem~D.  
In Theorem \ref{HLP-th4}, let $\lambda=0$, we can get
$$
M_f(r) +\left(\frac{1}{1+a_0}+\frac{r}{1-r}\right)\|f_0\|_r^2\leq 1 
$$
for $|z|=r\le r_{\lambda}=1/(5-2a_0)\in[1/5, 1/3)$ and $0.24683\in [1/5, 1/3)$. Moreover, in Theorem \ref{HLP-th4}, when $a_0\approx0.47431$,
we can get that (\ref{liu37}) $D_\lambda(r)\le1$ for $|z|=r\le r_*\approx0.24683$. And when $a_0$ takes other values, there exists $1/5\le r_{\lambda}<1/3$. Note also that
$$
\frac{1}{1+a_0}+\frac{r}{1-r}  \geq \frac{1}{2} \ \mbox{for all $r\geq0$ and $a\in [0,1)$}
$$
 and so Theorem \ref{HLP-th4} may be compared with \cite{KP2017} for $f \in \mathcal{B}.$
\er

\section{Bohr inequalities of alternating series $A_f(r), f \in \mathcal{B}.$}\label{LLP-sec4}
The definitions of alternating series $A_f(r)$, $A_{f_{0}}(r)$ and $A_{f_{m}}(r)$ are given in Section \ref{LSL-subsec1-2}.
The main results of this section are following: Firstly, estimate the upper bounds of the odd and even terms of the function
$f(z)=\sum_{n=0}^{\infty}a_{pn+m}z^{pn+m}\in\mathcal{B}$, then use $\{\zeta_n(r)\}_{n=0}^{\infty}$ to redefine
the bohr radius of the alternating series $A_f(r)$, and finally generalize the Theorem~E. 
Most of the conclusions are accurate.

\blem\label{HLP-th5}
Assume that $\zeta=\{\zeta_n(r)\}_{n=0}^{\infty}$ belonging to $\mathcal{F}$. If $f(z)=\sum_{n=0}^{\infty}a_{pn+m}z^{pn+m}\in\mathcal{B}$,
where $p\in \mathbb{N}$ and $0\le m \le p$, then

\begin{enumerate}
\item[{\rm (i)}] we have
\begin{eqnarray}
  &&\sum_{n=0}^{\infty}|a_{2np+m}|\zeta_{2n}(r^p)+\sum_{n=1}^{\infty}|a_{np+m}|^2\left(\frac{\zeta_{2n}(r^p)}{1+|a_m|}+\sum_{k=n}^{\infty}\zeta_{2(k+1)}(r^p)\right)\nonumber\\
  && \hspace{5cm}\le|a_m|\zeta_0(r^p)+(1-|a_m|^2)\sum_{n=1}^{\infty}\zeta_{2n}(r^p)
\label{liu41}
\end{eqnarray}
 for $r\in[0,1).$ The result is sharp.


\item[{\rm (ii)}] we have
\begin{equation}
  \sum_{n=1}^{\infty}|a_{(2n-1)p+m}|\zeta_{2n-1}(r^p)+\sum_{n=0}^{\infty}|a_{np+m}|^2\sum_{k=n}^{\infty}\zeta_{2k+1}(r^p)\le\sum_{n=1}^{\infty}\zeta_{2n-1}(r^p)
\label{liu42}
\end{equation}
 for $r\in[0,1).$ This result is sharp.
\end{enumerate}
\elem

\br
\begin{enumerate}
\item When $\zeta_{n}(r)=r^{n}\ (n\ge0)$, $m=0$ and $p=1$, the inequality \eqref{liu41} is the same as Remark \ref{HLP-re1}(i) with $k=2$.
\item When $\zeta_{n}(r)=r^{n}\ (n\ge1)$, $m=0$ and $p=1$, the inequality \eqref{liu42} is the same as Remark \ref{HLP-re1}(ii).
\end{enumerate}
\er

\bpf
(i) We may write $f$ as $f(z)=z^m\sum_{n=0}^{\infty}a_{pn+m}(z^{p})^n=z^mg(z^p)$, where $g \in\mathcal{B}$ and $g(s)=\sum_{n=0}^{\infty}b_ns^n$ with $a_{pn+m}=b_n$. 
If we apply Lemma~G 
 to the function $g(s)$
($s=z^p$), we get the desired inequality \eqref{liu41} with let $\zeta_{2n+1}(r)=0\ (n\ge0)$:
\begin{equation*}
\sum_{n=0}^{\infty}|b_{2n}|\zeta_{2n}(\rho)+\sum_{n=1}^{\infty}|b_{n}|^2\left(\frac{\zeta_{2n}(\rho)}{1+|b_0|}+\sum_{k=n}^{\infty}\zeta_{2(k+1)}(\rho)\right) \le|b_0|\zeta_0(\rho)+(1-|b_0|^2)\sum_{n=1}^{\infty}\zeta_{2n}(\rho).
\end{equation*}
where $\rho =r^p$.

To prove the sharpness, we let $a\in [0,1)$ and consider the function
\beq \label{1(a)}
f(z)=z^m\left(\frac{a-z^p}{1-az^p}\right)=az^m-(1-a^{2})\sum_{n=1}^{\infty}a^{n-1}z^{np+m},\quad  z\in \mathbb{D}.
\eeq
For this function, using the proof method as \cite[Theorem 1]{PVW2021}, it can be easily seen that inequality in \eqref{liu41} becomes an equality.
The proof of Lemma \ref{HLP-th5}(i) is completed.

(ii) In this case, we only need to apply Lemma~G 
to function $g(s)$, where $f(z)=z^m g(z^p)$ with $g(s)=\sum_{n=0}^{\infty}b_ns^n$  and
$b_n=a_{np+m}$, and let $\zeta_{2n}(r)=0\ (n\ge0)$. The rest of the process is exactly the same as Lemma \ref{HLP-th5}(i). 
\epf

 \bcor\label{HLP-cor4}
Assume that $\zeta=\{\zeta_n(r)\}_{n=0}^{\infty}$ belongs to $\mathcal{F}$, and $f \in \mathcal{B}$ with $f(z)=\sum_{n=0}^{\infty}a_{pn+m}z^{pn+m}$, where $p\in \mathbb{N}$ and $0\le m \le p$. Then
\begin{eqnarray}
 &&\sum_{n=1}^{\infty}|a_{np+m}|\zeta_{n}(r^p)+\sum_{n=1}^{\infty}|a_{np+m}|^2\left(\frac{\zeta_{2n}(r^p)}{1+|a_m|}+\sum_{k=2n+1}^{\infty}\zeta_{k}(r^p)\right)\nonumber\\ && \hspace{6cm}\le(1-|a_m|^2)\sum_{n=1}^{\infty}\zeta_{n}(r^p)
\label{liu43}
\end{eqnarray}
holds for $r\in[0,1)$. The result is sharp for $f(z)=z^m\frac{a-z^p}{1-az^p}$ with $a\in[0,1)$.
\ecor

\bpf
Add  (\ref{liu41}) and (\ref{liu42}), and then organize them to get (\ref{liu43}). To prove that the equality sign in (\ref{liu43}) holds, we let $a\in[0,1)$ and consider the function $f$ defined by \eqref{1(a)}, and follow the usual procedure.
\epf

\br
When $\zeta_{n}(r^p)=r^{np}\ (n\ge0)$ with $m=0$ and $p=1$, then \eqref{liu43} is the same as Remark \ref{HLP-re1}(iii).
\er

\br\label{HLP-re6}
Corollary \ref{HLP-cor4} establishes a connection between (\ref{liu41}), (\ref{liu42}) and (\ref{liu43}). As long as you know two of these formulas, you can deduce another formula. In particular, when $\zeta_{n}(r^p)=r^{np}\ (n\ge0)$, $p=1$ and $m=0$, the relationships between $M_f(r)$ and odd and even terms are established.
\er

In the following, Lemma \ref{HLP-th5} and Corollary \ref{HLP-cor4} will be used to explore the redefinition of the Bohr radius of alternating series $A_f(r)$.

\subsection{Basic notation}
We first introduce some  definitions and notations: For all $n\ge0$, let $\zeta_n^*(r)=C_n r^{\tau_n},$ where $C_n=C(n)\geq 0,\,$ 
and $\{\tau_n\}_{n \geq 1}$ is a strictly increasing sequence of natural numbers,
then $\{\zeta^*_n(r)\}_{n=0}^{\infty}\subset\{\zeta_n(r)\}_{n=0}^{\infty}\in\mathcal{F}$. For example, if $f(z)=z^m\sum_{n=0}^{\infty}a_{np+m}z^{np}$ for some $m\geq 0$ and $p\geq 1$,
then we may introduce
\begin{equation*}
  A^*_{f_{m}}(r)=\zeta_m^*(r)\sum_{n=1}^{\infty}(-1)^{\tau_{np+m}}|a_{np+m}|\zeta_n^*(r^p),
\end{equation*}
where  $\zeta_m^*(r)\zeta_n^*(r^p)=C_{np+m} r^{\tau_{np+m}}$ with $C_{np+m} =C(n,p,m)\ge0$, and $\zeta_0^*(r)=1$ for all $r\in[0,1)$. Moreover, if we let $\zeta_m^*(r)=r^m$ then $\zeta_n^*(r^p)=r^{np}$ so that, with $\tau_{np+m}=np+m,$ $A^*_{f_{m}}(r)$ reduces to
\begin{equation*}
A^*_{f_{m}}(r)=\sum_{n=1}^{\infty}(-1)^{np+m}|a_{np+m}|r^{np+m}=:A_{f_{m}}(r).
\end{equation*}

\bthm\label{HLP-th6}
Assume that $\zeta^*=\{\zeta^*_n(r)\}_{n=0}^{\infty}$ belonging to $\mathcal{F}$ is a decreasing sequence of nonnegative functions, $f \in \mathcal{B}$ and $f(z)=z^m\sum_{n=0}^{\infty}a_{np+m}(z^{p})^n$, where $p\in \mathbb{N}$, $0\le m \le p$. For convenience, we let
\begin{equation*}
 A^*_{f_{m}}(r) =\zeta^*_m(r) B^*_{f_{m}}(r),\quad B^*_{f_{m}}(r)=\sum_{n=1}^{\infty}(-1)^{\tau_{np+m}}|a_{np+m}|\zeta^*_n(r^p),
\end{equation*}
\begin{equation*}
 C^*_{f_{m}}(r)= \zeta^*_m(r) \bigg[B^*_{f_{m}}(r) +(-1)^{\tau_{(2n-1)p+m}} \sum_{n=0}^{\infty}|a_{np+m}|^2\sum_{k=n}^{\infty}\zeta^*_{2k+1}(r^p) \bigg],
 ~\mbox{ and }
\end{equation*}
\begin{equation*}
 D^*_{f_{m}}(r)= \zeta^*_m(r) \bigg[B^*_{f_{m}}(r) +(-1)^{\tau_{np+m}} \sum_{n=1}^{\infty}|a_{np+m}|^2\left(\frac{\zeta^*_{2n}(r^p)}{1+|a_m|}+\sum_{k=2n+1}^{\infty}\zeta^*_{k}(r^p)\right)\bigg].
\end{equation*}
We have the following:
\begin{itemize}
\item[(I)] If the degree $\tau_{2np+m}$ of $\zeta^*_{2n}(r^p)\zeta^*_m(r)$ is even (resp. odd) and the degree $\tau_{(2n-1)p+m}$ of $\zeta^*_{2n-1}(r^p)\zeta^*_m(r)$ is odd (resp. even), 
then   the sharp inequality
\begin{equation}
 \bigg|C^* _{f_{m}}(r) \bigg| \le 1 ~\mbox{ for $|z|=r\le r_*$}
\label{liu44}
\end{equation}
holds, where $r_{*}$ is the root of the equation $\zeta^*_m(r)\sum_{n=1}^{\infty}\zeta^*_{2n-1}(r^p)=1$.
\item[(II)] If the degrees of $\zeta^*_{2n}(r^p)\zeta^*_m(r)$ and $\zeta^*_{2n-1}(r^p)\zeta^*_m(r)$ are both even
(resp.  both odd), then the sharp inequality
\begin{equation}
 \bigg| D^*_{f_{m}}(r) \bigg| \le 1 ~\mbox{ for $|z|=r\le R_*$}
\label{liu45}
\end{equation}
holds,  where $R_{*}$ is the root of the equation $\zeta^*_m(r)\sum_{n=1}^{\infty}\zeta^*_{n}(r^p)=1$.
\end{itemize}
\ethm

\bpf
(I) Assume first that the degree of $\zeta^*_{2n}(r^p)\zeta^*_m(r)$ is even and the degree of $\zeta^*_{2n-1}(r^p)\zeta^*_m(r)$ is odd, both for all $n\geq 1$. Then, we may write
$$
A^*_{f_{m}}(r)  =\zeta^*_m(r)B^*_{f_{m}}(r) =\zeta^*_m(r)\bigg[\sum_{n=1}^{\infty}|a_{2np+m}|\zeta^*_{2n}(r^p)-\sum_{n=1}^{\infty}|a_{(2n-1)p+m}|\zeta^*_{2n-1}(r^p)\bigg].
$$
Because all $\zeta_n^*\in\mathcal{F}~ (n\ge0)$, by Lemma \ref{HLP-th5}(ii) and the definition of  $B^*_{f_{m}}(r),$ we have
\begin{align}
B^*_{f_{m}}(r) &-\sum_{n=0}^{\infty}|a_{np+m}|^2\sum_{k=n}^{\infty}\zeta^*_{2k+1}(r^p) \nonumber  \\
&\ge-\sum_{n=1}^{\infty}|a_{(2n-1)p+m}|\zeta^*_{2n-1}(r^p)-\sum_{n=0}^{\infty}|a_{np+m}|^2\sum_{k=n}^{\infty}\zeta^*_{2k+1}(r^p) \nonumber  \\
& \ge -\sum_{n=1}^{\infty}\zeta^*_{2n-1}(r^p). \label{Liu4.7}
\end{align}
Again by Lemma  \ref{HLP-th5}(i), since $\{\zeta^*_n(r)\}_{n=0}^{\infty}$ is a decreasing sequence, we have $\zeta^*_{2n}(r^p)\le\zeta^*_{2n-1}(r^p)$ for all $r\in[0,1)$ and thus,
\begin{align}
B^*_{f_{m}}(r) &-\sum_{n=0}^{\infty}|a_{np+m}|^2\sum_{k=n}^{\infty}\zeta^*_{2k+1}(r^p)\le \sum_{n=1}^{\infty}|a_{2np+m}|\zeta^*_{2n}(r^p) \nonumber\\
&\le(1-|a_m|^2)\sum_{n=1}^{\infty}\zeta^*_{2n}(r^p)-\sum_{n=1}^{\infty}|a_{np+m}|^2\bigg(\frac{\zeta^*_{2n}(r^p)}{1+|a_m|}+\sum_{k=n}^{\infty}\zeta^*_{2(k+1)}(r^p)\bigg) \nonumber\\
&\le \sum_{n=1}^{\infty}\zeta^*_{2n}(r^p)\le\sum_{n=1}^{\infty}\zeta^*_{2n-1}(r^p). \nonumber
\end{align}
which by, combining with the inequality \eqref{Liu4.7} gives
\beqq
\bigg| C^*_{f_{m}}(r) \bigg| \le \zeta^*_m(r)\sum_{n=1}^{\infty}\zeta^*_{2n-1}(r^p)=h_1(r)+1,
\eeqq
which is less than or equal to $1$,   provided $h_1(r)\leq 0$, where
$$h_1(r):=\zeta^*_m(r)\sum_{n=1}^{\infty}\zeta^*_{2n-1}(r^p)-1.
$$
Note that 
$h_1(r)=\sum_{n=1}^{\infty}C((2n-1),p,m)r^{\tau_{(2n-1)p+m}}-1$ and $h_1'(r)>0$ 
for $r\in[0,1].$ Since the degree of $\zeta^*_{2n-1}(r^p)\zeta^*_m(r)$ is odd, it follows that
$$
h_1(0)=\zeta^*_m(0)\sum_{n=1}^{\infty}\zeta^*_{2n-1}(0)-1=-1  <0.
$$
According to the intermediate value theorem, we deduce that $h_1(r)\leq 0$  for $|z|=r\leq r_*$.

To prove the sharpness of the first case, we consider the function $f(z)=z^{p+m}=z^mz^p$ (that is just $|a_{p+m}|=1$ and $|a_{np+m}|=0$, $n\neq1$), because the degree $\tau_{(2n-1)p+m}$ of $\zeta^*_{2n-1}(r^p)\zeta^*_m(r)$ is odd, we see that $A^*_{f_{m}}(r)=-\zeta^*_1(r^p)\zeta^*_m(r)$. In this case we get
\begin{equation*}
\bigg| C^*_{f_{m}}(r) \bigg|=\bigg|-\zeta^*_1(r^p)\zeta^*_m(r)-\zeta^*_m(r)\sum_{k=1}^{\infty}\zeta^*_{2k+1}(r^p)\bigg|=\zeta^*_m(r)\sum_{n=1}^{\infty}\zeta^*_{2n-1}(r^p).
\end{equation*}
which is bigger than 1 for $|z|=r > r_*$. This completes the proof of the sharpness for the first case.

Next we may suppose that the degree  $\tau_{2np+m}$ of $\zeta^*_{2n}(r^p)\zeta^*_m(r)$ is odd and the degree $\tau_{(2n-1)p+m}$ of $\zeta^*_{2n-1}(r^p)\zeta^*_m(r)$ is even for all $n\ge0$. In this case, we see that
\begin{equation*}
A^*_{f_{m}}(r)=\zeta^*_m(r)\bigg[-\sum_{n=1}^{\infty}|a_{2np+m}|\zeta^*_{2n}(r^p)+\sum_{n=1}^{\infty}|a_{(2n-1)p+m}|\zeta^*_{2n-1}(r^p)\bigg].
\end{equation*}
Next combined with Lemma \ref{HLP-th5}, the rest of the proof is similar to the first part of Case I above, except that the corresponding formula
has the opposite sign. To prove the sharpness, we also consider the function $f(z)=z^{m+p}$.

This completes the proof of Case I.

(II) Assume that the degrees of $\zeta^*_{2n}(r^p)\zeta^*_m(r)$ and $\zeta^*_{2n-1}(r^p)\zeta^*_m(r)$ are both even for all $n\ge0$. This means that all $\tau_{np+m}$ are even.
In this case, we see that
$$
A^*_{f_{m}}(r) =\zeta^*_m(r)B^*_{f_{m}}(r) =\zeta^*_m(r)\sum_{n=1}^{\infty}|a_{np+m}|\zeta^*_{n}(r^p).
$$
Thus, since all $\zeta_n^*\in\mathcal{F}~ (n\ge0),$  Corollary \ref{HLP-cor4} yields that
\begin{align}
 \left | D^*_{f_{m}}(r) \right |
&=\zeta^*_m(r)\bigg[\sum_{n=1}^{\infty}|a_{np+m}|\zeta^*_{n}(r^p)+\sum_{n=1}^{\infty}|a_{np+m}|^2\bigg(\frac{\zeta^*_{2n}(r^p)}{1+|a_m|}+\sum_{k=2n+1}^{\infty}\zeta^*_{k}(r^p)\bigg)\bigg]
\nonumber \\
&\le \zeta^*_m(r)(1-|a_m|^2)\sum_{n=1}^{\infty}\zeta^*_{n}(r^p)\le\zeta^*_m(r)\sum_{n=1}^{\infty}\zeta^*_{n}(r^p)=h_2(r)+1
\label{liu48}
\end{align}
which is less than or equal to $1$ provided $h_2(r) \leq 0,$ where $h_2(r):=\zeta^*_m(r)\sum_{n=1}^{\infty}\zeta^*_{n}(r^p)-1.$ $h_2(r)$  is increasing on $[0,1].$ Since each $\tau_{np+m}$ is different, $\zeta^*_0(r)=1$ and $\{\zeta^*_n(r)\}_{n=0}^{\infty}$ is a decreasing sequence   from the hypothesis, we find that $h_2(0)\le\zeta^*_m(0)\zeta^*_{j}(0)-1\le\zeta^*_0(0)\zeta^*_{0}(0)-1=0$. 
By intermediate value Theorem, we have  $h_2(r)\le0$  whenever $r\le R_*$.

To prove the sharpness, we consider the function $f(z)=z^{p+m}=z^mz^p$ and note that the degree of $\zeta^*_{2n-1}(r^p)\zeta^*_m(r)$ is even number, so that $A^*_{f_{m}}(r)=\zeta^*_1(r^p)\zeta^*_m(r)$. In this case, it follows easily that
\begin{equation}\nonumber
\begin{split}
 \left | D^*_{f_{m}}(r) \right |
=\zeta^*_m(r)\bigg[\zeta^*_1(r^p)+\zeta^*_2(r^p)+\sum_{n=3}^{\infty}\zeta^*_{k}(r^p)\bigg]=\zeta^*_m(r)\sum_{n=1}^{\infty}\zeta^*_{n}(r^p)
\end{split}
\end{equation}
which is bigger than $1$ for $|z|=r > R_*.$   This completes the proof in this case.

In the second part of Case II, we let   the degrees of $\zeta^*_{2n}(r^p)\zeta^*_m(r)$ and $\zeta^*_{2n-1}(r^p)\zeta^*_m(r)$  are both odd for all $n\ge0$. Then, we write
\begin{equation*}
A^*_{f_{m}}(r)=\zeta^*_m(r)\bigg[-\sum_{n=1}^{\infty}|a_{2np+m}|\zeta^*_{2n}(r^p)-\sum_{n=1}^{\infty}|a_{(2n-1)p+m}|\zeta^*_{2n-1}(r^p)\bigg].
\end{equation*}
Next combined with Corollary \ref{HLP-cor4}, the rest of proof is similar to the first part of Case II, except that the corresponding formula has the opposite sign. To prove the sharpness, we just consider again the function $f(z)=z^{m+p}$.
\epf

\begin{example} \label{HLP-ex7} Using Theorem \ref{HLP-th6} and letting $\zeta^*_m(r)=r^m$ and $\zeta^*_n(r^p)=r^{np} \ (n\ge0)$, we deduce the following  three examples
whenever $f \in\mathcal{B}.$ 

\begin{enumerate}
\item[(1)] Set $p=1$ and $m=0$ so that $f$ has the form $f(z)=\sum_{n=0}^{\infty}a_nz^n$ and that $A_{f_0}(z)=\sum_{n=1}^{\infty}(-1)^n a_nz^n.$ By the first case of Theorem \ref{HLP-th6}(I), we get that
$$\left | A_{f_0}(r)-\frac{r}{1-r^{2}}\sum_{n=0}^{\infty}|a_{n}|^2r^{2n}\right | = \left | A_{f_0}(r)-\frac{r}{1-r^{2}}\Vert f \Vert_r ^2 \right | \le 1 \ \mbox{ for $\ds |z|=r\le \frac{\sqrt{5}-1}{2}$.}
$$
This result is sharp for $f(z)=z$. The earlier known result is $|A_{f_{0}}|\le1$ holds for $|z|=r\le\frac{\sqrt{5}-1}{2}$.

\item[(2)] Set $p=k$ is even and $m=0$ so that $f(z)=\sum_{n=0}^{\infty}a_{kn}z^{kn}$ is even function. By the first case of Theorem \ref{HLP-th6}(II), we get that
$$\left | \sum_{n=1}^{\infty}|a_{kn}|r^{kn}+\left(\frac{1}{1+|a_0|}+\frac{r^k}{1-r^k}\right)\sum_{n=1}^{\infty}|a_{kn}|^2r^{2kn}\right |\le 1
$$
for $|z|=r\le\sqrt[k]{\frac{1}{2}}$. The result is sharp for $f(z)=z^k$.

\item[(3)] Finally,   set $p=1$, $m=0$ so that $f(z)=\sum_{n=0}^{\infty}a_{n}z^{n}$, but now let $\zeta^*_n(r)=\frac{r^n}{n+1} ~(n\ge0)$, then by Theorem \ref{HLP-th6}(I), we have
$$\left | \sum_{n=1}^{\infty}(-1)^{n+1}\frac{|a_n|}{n+1}r^{n}+\sum_{n=1}^{\infty}|a_n|^2\frac{1}{r}\int_{0}^{r}\frac{t^{2n+1}}{1-t^2}dt\right |\le 1.
$$
holds for $|z|=r\le R$, which is the root of $-\log(1-r^2)=2r$. The result is sharp for $f(z)=z$.
\end{enumerate}
\end{example}

By replacing $A^*_{f_{m}}(r)$ in (\ref{liu45}) of Theorem \ref{HLP-th6}(II) with $A^*_{f}(r)=(-1)^{\tau_{m}}\zeta^*_m(r)|a_m|+A^*_{f_{m}}(r)$, we get a generalization of Theorem \ref{HLP-th6}(II) with the first coefficient.

\bcor\label{HLP-cor4b}
  If the degrees of $\zeta^*_{2n}(r^p)\zeta^*_m(r)$ and $\zeta^*_{2n-1}(r^p)\zeta^*_m(r)$ are both even
(resp. both odd), we have
\begin{equation}
 \bigg|A^*_{f}(r)+(-1)^{\tau_{np+m}}\zeta^*_m(r)\sum_{n=1}^{\infty}|a_{np+m}|^2\bigg(\frac{\zeta^*_{2n}(r^p)}{1+|a_m|}+\sum_{k=2n+1}^{\infty}\zeta^*_k(r^p)\bigg)\bigg| \le1
\label{liu49a}
\end{equation}
holds for $|z|=r\le \tilde{r}$, which is the root of $\zeta^*_m(r)+4\sum_{n=1}^{\infty}\zeta^*_n(r^p)[\sum_{n=1}^{\infty}\zeta^*_n(r^p)\zeta^*_m(r)-1]=0$. This result is sharp for $f(z)=z^m\frac{z^p-a}{1-az^p}$ with $0<a<1$.
\ecor

\bpf
 With the help of Corollary \ref{HLP-cor4}, the proof is similar as in \cite[Lemma 2.2]{ABS2017}.
\epf

\begin{example} \label{HLP-ex8} Let $\zeta^*_n(r^p)=r^{np} ~(n\ge0)$ and $\zeta^*_m(r)=r^m$ in Corollary \ref{HLP-cor4b}. Then we obtain the
following two examples, which can be regarded as the promotion of \cite[Lemma 2.1, Lemma 2.2]{ABS2017}.

\begin{enumerate}
\item[(i)] Set $m=0$ and $p=2$, so that $f(z)=\sum_{n=0}^{\infty}a_{2n}z^{2n}\in\mathcal{B}$. By Corollary \ref{HLP-cor4b}, we have
\begin{equation*}
 \sum_{n=0}^{\infty}|a_{2n}|r^{2n}+\bigg(\frac{1}{1+|a_0|}+\frac{r^2}{1-r^2}\bigg)\sum_{n=1}^{\infty}|a_{2n}|^2r^{4n}\le1 ~\mbox{ for $\ds |z|=r\le \frac{1}{\sqrt{3}}$.}
\end{equation*}
This result is sharp for $f(z)=\frac{z^2-a}{1-az^2}$ with $0<a<1$.

\item[(ii)] Set $m=1$ and $p=2$, so $f(z)=\sum_{n=0}^{\infty}a_{2n+1}z^{2n+1}\in\mathcal{B}$. By Corollary \ref{HLP-cor4b}, we have
\begin{equation*}
 \sum_{n=0}^{\infty}|a_{2n+1}|r^{2n+1}+\bigg(\frac{r}{1+|a_1|}+\frac{r^3}{1-r^2}\bigg)\sum_{n=1}^{\infty}|a_{2n+1}|^2r^{4n}\le1 ~\mbox{$|z|=r\le \tilde{r_2} $,}
\end{equation*}
where $\tilde{r_2}\approx 0.731348$ is the root of $5r^4+4r^3-2r^2-4r+1=0$. This result is sharp for $f(z)=z\frac{z^2-a}{1-az^2}$ with $0<a<1$.
\end{enumerate}
\end{example}

Finally, we present an improved version of Theorem~E. 
In this case, we only need to discuss the case of $\zeta_n(r)=r^{n}$. For convenience, let
\begin{equation*}
E_{f_{m}}(r)= r^m \bigg[\sum_{n=1}^{\infty}(-1)^{np+m}|a_{np+m}|r^{np} +(-1)^m  \left(\frac{1}{1+|a_m|}+\frac{r^{2p}}{1-r^{2p}}\right)\sum_{n=1}^{\infty}|a_{np+m}|^2r^{2np} \bigg].
\end{equation*}

\bthm\label{HLP-th7}
Suppose that $f \in \mathcal{B}$ and $f(z)=z^m\sum_{n=0}^{\infty}a_{np+m}z^{np}$, where $0\le m\le p$, $m \in \mathbb{N}$ is even and $p \in \mathbb{N}$ is odd. Then
\begin{equation}
 \left | \ |f(z)|+  E_{f_{m}}(r) \ \right | \le 1
\label{liu49}
\end{equation}
 for $|z|=r\le r_{p,m}$, where $r_{p,m}$ is the minimal positive root of $r^{3p-m}-2r^{3p}-3r^{2p}-r^{p-m}+1=0$, and $E_{f_{m}}(r)$ is defined as above.
\ethm

\bpf
By assumption, each $2np+m$ ($n\ge0$) is even and each $(2n-1)p+m$ ($n\ge1$) is odd, as $m$ is even and $p$ is odd. Thus, we write
$$A_{f_{m}}(r)= r^m \sum_{n=1}^{\infty}(-1)^{np+m}|a_{np+m}|r^{np}=r^m\bigg[\sum_{n=1}^{\infty}|a_{2np+m}|r^{2np}-\sum_{n=1}^{\infty}|a_{(2n-1)p+m}|r^{(2n-1)p}\bigg].
$$
 Next  by Lemma \ref{HLP-th5}(i) with $\zeta_{2n}(r^p)=r^{2np} (n\ge0)$ and $m$ is even, then
\beqq
E_{f_{m}}(r) &=&r^m\bigg[\sum_{n=1}^{\infty}(-1)^{np+m}|a_{np+m}|r^{np}+\left(\frac{1}{1+|a_m|}+\frac{r^{2p}}{1-r^{2p}}\right)\sum_{n=1}^{\infty}|a_{np+m}|^2r^{2np}\bigg]\nonumber\\
 &\le&r^m\bigg[\sum_{n=1}^{\infty}|a_{2np+m}|r^{2np}+ \left(\frac{1}{1+|a_m|}+\frac{r^{2p}}{1-r^{2p}}\right)\sum_{n=1}^{\infty}|a_{np+m}|^2r^{2np}\bigg]\nonumber\\
 &\le& (1-|a_m|^2)\frac{r^{2p+m}}{1-r^{2p}},\nonumber\\
\eeqq
so that
\be\label{liu410}
|f(z)|+  E_{f_{m}}(r)  \leq r^m\frac{r^p+|a_m|}{1+r^p|a_m|}+(1-|a_m|^2)\frac{r^{2p+m}}{1-r^{2p}}.
\ee
To find a lower bound for $|f(z)|+E_{f_{m}}(r)$,  by Lemma \ref{HLP-th5}(ii) with $\zeta_{2n-1}(r^p)=r^{(2n-1)p}~ (n\ge1)$, we obtain
\beqq
|f(z)|+  E_{f_{m}}(r) &\geq& -r^m\sum_{n=1}^{\infty}|a_{(2n-1)p+m}|r^{(2n-1)p}\\
&\ge & -\frac{r^{p+m}}{1-r^{2p}}+\frac{r^{p+m}}{1-r^{2p}}\sum_{n=0}^{\infty}|a_{np+m}|^2r^{2np}\\
&\ge& -\frac{r^{p+m}}{1-r^{2p}}(1-|a_m|^2)= -\frac{r^{2p+m}}{1-r^{2p}}(1-|a_m|^2)-r^m   \cdot\frac{r^p(1-|a_m|^2)}{1+r^p},
\eeqq
  where the third inequality above, we have used the identity (for $x=r^p$)
$$\frac{1}{1-x^2}=\frac{x}{1-x^2}+\frac{1}{1+x} \ \mbox{for $x \in [0,1)$},
$$
and because
$$\frac{r^p+|a_m|}{1+r^p|a_m|}-\frac{r^p(1-|a_m|^2)}{1+r^p}=\frac{(1-|a_m|+|a_m|^3)r^{2p}+|a_m|[(1+|a_m|)r^p+1]}{(1+r^p|a_m|)(1+r^p)}\geq 0,
$$
we get that
\be\label{liu410a}
|f(z)|+E_{f_{m}}(r)\geq -\frac{r^{2p+m}}{1-r^{2p}}(1-|a_m|^2)-r^m\frac{r^p+|a_m|}{1+r^p|a_m|},
\ee
and thus, combining \eqref{liu410} and \eqref{liu410a} gives
\be\label{liu410c}
\left|  |f(z)|+ E_{f_{m}}(r)  \right | \leq r^m\frac{r^p+|a_m|}{1+r^p|a_m|}+(1-|a_m|^2)\frac{r^{2p+m}}{1-r^{2p}}.
\ee
According to Lemma \ref{HLP-lem5}, we find that the right hand side of (\ref{liu410c}) is smaller than or equal to $1$ for $|z|=r\le r_{p,m}$, where $r_{p,m}$ is as in the statement. This completes the proof of inequality (\ref{liu49}).
\epf

\bcor\label{HLP-cor5}According to Theorem \ref{HLP-th7}, we can present an improved version of Theorem~E. 
Suppose that $f\in\mathcal{B}$ and $f(z)=\sum_{n=0}^{\infty}a_nz^n$. Then
\begin{equation*}
 \left | |f(z)|+A_{f_{0}}(r)+\left(\frac{1}{1+|a_0|}+\frac{r^2}{1-r^2}\right)\sum_{n=1}^{\infty}|a_{n}|^2r^{2n}\right | \le 1
\end{equation*}
holds for $|z|=r\le\sqrt{2}-1.$
\ecor

It is worth pointing out that we have not proved that the number $r=\sqrt{2}-1$ is the best possible in Theorem \ref{HLP-th7} and therefore,
the following problem remains open.

\bprob\label{HLP-prob2} 
Find the largest radius $r_0$ for the class $\mathcal{B}$ of analytic functions $f(z)=\sum_{n=0}^{\infty}a_nz^n$ in $\mathbb{D}$  such that
\begin{equation*}
\bigg||f(z)|+\sum_{n=1}^{\infty}(-1)^n|a_n|r^n+\left(\frac{1}{1+|a_0|}+\frac{r^2}{1-r^2}\right)\Vert f_0 \Vert _{r}^{2}\bigg|\leq 1~~\mbox{ for }~~|z|=r\leq r_0.
\end{equation*}
\eprob

\bcor\label{HLP-cor6}
Suppose that $f(z)\in\mathcal{B}$ and $f(z)=\sum_{n=0}^{\infty}a_{np+m}z^{np+m}$, where $m \in \mathbb{N}$ is even and $p \in \mathbb{N}$ is odd. Then
\begin{equation*}
 \left | \frac{r^{m+p}}{1+r^p}+A_{f_{m}}(r)+r^m\left(\frac{1}{1+|a_m|}+\frac{r^{2p}}{1-r^{2p}}\right)\sum_{n=1}^{\infty}|a_{np+m}|^2r^{2np}\right | \le 1
\label{liu52}
\end{equation*}
for $|z|=r\le R_{p,m}$, where $R_{p,m}$ is the root of $r^{2p}+r^{p+m}-1=0$. The result is sharp for $f(z)=z^{2p+m}$.
\ecor

\bpf
Although Theorem \ref{HLP-th7} cannot provide an accurate result, the result becomes accurate if we replace the term $|f(z)|$ by $\frac{r^{m+p}}{1+r^p}$. The method of proof of Corollary \ref{HLP-cor6} is exactly the same as Theorem \ref{HLP-th7}, so it will not be repeated. Note that condition $\frac{r^{m+p}}{1+r^p}$ in Corollary \ref{HLP-cor6} is necessary, otherwise the scaling of the lower bound will not hold.
\epf


\subsection*{Acknowledgments}
The work of the first two authors are supported by Guangdong Natural Science Foundations (Grant No. 2021A1515010058). The work of the third author was supported by Mathematical Research Impact Centric Support (MATRICS) of the Department of Science and Technology (DST), India  (MTR/2017/000367).

\subsection*{Conflict of Interests}
The authors declare that they have no conflict of interest, regarding the publication of this paper.
\subsection*{Data Availability Statement}
The authors declare that this research is purely theoretical and does not associate with any data.

\end{document}